\documentclass{article}
\usepackage[T2A]{fontenc}
\usepackage[utf8]{inputenc}
\usepackage[english,russian]{babel}
\usepackage[tbtags]{amsmath}
\usepackage{amsfonts,amssymb,mathrsfs,amscd}
\usepackage[hyper]{msb-a}

\makeatletter
\gdef\No{{\select@language{russian}\textnumero}}
\makeatother

\JournalName{}
\numberwithin{equation}{section}
\def \Z {{\mathbf {Z}}}

\def \B {{\cal B}}
\def \T {{\mathbf  T}}

\def\eps{\varepsilon}

\usepackage{graphicx}
\graphicspath{{pictures/}}
\DeclareGraphicsExtensions{.pdf,.png,.jpg}
\begin{document}

\title{Кратное перемешивание в эргодической теории} 
%Cемидесят пять лет  проблеме В.А. Рохлина}
\author[V.\,V.~Ryzhikov]{В.\,В.~Рыжиков}
\address{Московский государственный университет}
\email{vryzh@mail.ru}

\date{2024}
\udk{517.987}

\maketitle

\begin{fulltext}

\begin{abstract}{

In 1949 V.A. Rokhlin introduced invariants into ergodic theory called k-fold mixing  and puzzled the mathematical community with the problem of the mismatch of these invariants.  Here's what Roslin wrote:
\it The proposed work arose from the author's attempts to solve the well-known spectral problem of the theory of dynamical systems: are there any metrically different dynamical systems with the same continuous (in particular, Lebesgue) spectrum?  Using new metric invariants introduced for this purpose, the author tried to find among the ergodic automorphisms of compact commutative groups automorphisms of various metric types. It turned out, however, that for all the indicated automorphisms the new invariants are exactly the same.\rm

We discuss    this long time open problem without  pretending to cover the topic.

В  1949 г. В.А.  Рохлин ввел в эргодическую теорию инварианты, названные  перемешиванием кратности $k$ и  озадачил математическое  сообщество  проблемой несовпадения этих инвариантов.   Мы напомним читателю некоторые результаты и методы их получения в связи с этой проблемой, которая 75 лет остается нерешенной.}
%Библиография: 29 названий.}
\end{abstract}
%\begin{keywords}
%     \end{keywords}
\markright{Кратное перемешивание в эргодической теории}

%%%%%%%%%%%%%%%%%%%%%%%%%%%%%%%%%%%%%%%%%%%%%%%%%%%%%%%%%%%%%%%%%%%%%%%%%%%
%\large
\section{Введение. Проблема Рохлина. Некоторые результаты. } 
В начале 1948 г. В.А.  Рохлин написал  
статью об эндоморфизмах компактных коммутативных групп,
которая была опубликовна в Известиях АН СССР  в 1949 году. Цитируем Рохлина \cite{Ro}:

\vspace{2mm}
\it Предлагаемая работа возникла из попыток автора решить известную
спектральную проблему теории динамических систем: существуют ли
метрически различные динамические системы с одним и тем же непрерывным 
(в частности, лебеговским) спектром? С помощью введенных
для этой цели новых метрических инвариантов  автор  пытался найти 
среди эргодических автоморфизмов компактных коммутативных групп 
автоморфизмы различных метрических типов. Оказалось, однако, что у всех указанных автоморфизмов новые инварианты совершенно одинаковы.\rm

\vspace{2mm}
Новые инварианты, предложенные в \cite{Ro}, называются свойствами 
перемешивания кратности $k\geq 1$.  Автоморфизм $T$ вероятностного пространства $(X,\B,\mu)$ обладает $k$-кратным перемешиванием, если для всяких измеримых множеств $A_0,A_1,\dots, A_k$ при $m_1,\dots m_k\, \to\,+\infty$  выполняется
$$\mu(A_0\cap T^{m_1}A_1\cap T^{m_1+m_2}A_2\cap \dots  T^{m_1+\dots+m_k}A_k)\ \to \ \mu(A_0)\mu(A_1)\dots \mu(A_k).$$
Обозначим это свойство и класс таких автоморфизмов через $Mix(k)$.

В \cite{Ro}  доказано кратное перемешивание  для   эргодических автоморфизмов компактных коммутативных групп (они сохраняют вероятностную меру Хаара). Хотя проблема Рохлина остается открытой, за прошедшие годы получен ряд  частных результатов, обнаружилась связь кратного перемешивания с различными инвариантами динамических систем. Проблема   стимулировала исследования самоприсоединений и марковских  сплетений тензорных произведений динамических систем.  Настоящая статья не является обзором, мы обсуждаем   некоторые методы доказательства кратного перемешивания.

%\section{Некоторые результаты}
Проблема Рохлина естественным образом обобщается на групповые действия.  Элегантное решение для $\Z^2$-действий, порожденное автоморфизмами компактной коммутативной  группы,   дал Ф. Ледрапье \cite{L}.  
На рисунке изображен типичный элемент  группы Ледрапье 
$H\subset \Z_2^{\Z^2}$.  Черным и красным квадратным точкам отвечает значение 0, а белым  -- значение 1.   $H$ состоит из гармонических последовательностей: значение в каждой точке равно $mod 2$ сумме значений в четырех соседних с ней точек.  Нити, которые видт читатель, автор обнаружил в 80-х годах, когда появились персональные компьютеры. Вероятно, в типичном случае (для почти всех элементов группы $H$)  одна бесконечная нить обеспечивает перколяцию, но это открытая проблема.

\begin{center}
\includegraphics{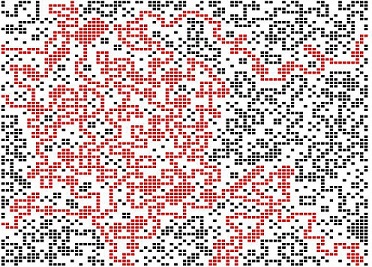}
\end{center}

 \bf Действие Ледрапье \rm порождено двумя автоморфизмами группы $H$:
вертикальным и горизонтальным сдвигами. Оно обладает перемешиванием кратности 3, но не обладает перемешиванием кратности 4. Интересное наблюдение появилось  в работе \cite{Ti}: максимальное отклонение от кратного перемешивания является причиной тривиальности централизатора действия.   Идея  Рохлина рассматривать для контрпримеров  действия,
порожденные автоморфизмами групп,   оказалась успешной в случае    $\Z^n$-действий для всех $n\neq 1$.    

Введеное Рохлиным  понятие  кратного перемешивания  обозначило направление от обычного перемешивания  в сторону более сильных перемешивающих инвариантов,  вершиной которых можно назвать K-перемешивание.  А.Н. Колмогоров  в 1958 г. решил    упомянутую Рохлиным спектральную проблему теории динамических систем (см. \cite{K}, \cite{KSF}): энтропия на единицу времени   могла  различить    бернуллиевские сдвиги, которые имеют одинаковый лебеговский спектр.  

\bf Асимметрия кратных асимптотические свойств. \rm Отметим, что первый пример пары неизоморфных автоморфизмов с одинаковым спектром был предложен Анзаи \cite{An}.  Искомой  парой являются  подходящие  косые произведения   $T$ и $T^{-1}$.  Что касается   идеи   Рохлина  использовать  кратные асимптотические свойства для различения автоморфизмов с одинаковым  спектром, в модифицированном  виде этот подход   был реализован (см. \cite{R23}): \it для типичного автоморфизма $T$  найдется последовательность $m_j\to\infty$ такая, что  для любого измеримого множества $A$ имеют место сходимости
$$4\mu(A\cap T^{m_j}A\cap T^{3m_j}A)\to \  \mu(A)+ \mu(A)^2 + 2\mu(A)^3,$$
$$4\mu(A\cap T^{-m_j}A\cap T^{-3m_j}A)\to \ 4\mu(A)^2.$$
\rm   Автоморфизмы $T$ и $T^{-1}$ имеют одинаковый спектр, но  разные метрические инвариантны. Указанные пределы не меняются при замене автоморфизма $T$ на сопряженый $S^{-1}TS$, поэтому 
$T$ и $T^{-1}$ не изоморфны.

\bf Ослабленная проблема Рохлина. \rm Задачу о кратном перемешивании можно ставить в классе слабо перемешивающих автоморфизмов.
Пусть для автоморфизма $T$ и последовательностей $m_i,n_i\to +\infty$ для  всех  $A,B\in \B$ выполнено
$$\mu(A\cap T^{m_i}B),\, \mu(A\cap T^{n_i}B),\, \mu(A\cap T^{n_i+m_i}B)\ \to\ \mu(A)\mu(B).$$

\it Будет ли для всех $A,B,C\in \B$ выполняется
$$\mu(A\cap T^{m_i} B\cap T^{m_i+n_i}C)\ \to \ \mu(A)\mu(B)\mu(C)?$$ 
\rm
Но и в этой ослабленной формулировке  проблема не становится проще.

\bf Теорема Оста. \rm Общий  результат  получил Б. Ост \cite{Ho}:
 \it  (слабо) перемешивающие автоморфизмы с сингулярным спектром (слабо) перемешивают  кратно. \rm

\vspace{2mm}
\bf  Естественность проблемы. \rm   Свойство  перемешивания  можно трактовать как  асимптотическую независимость события сейчас и  события в далеком прошлом.  Будет ли событие сейчас при выполнении указанного свойства и, что важно, при   стационарности случайного процесса    асимптотичеси не зависеть от пары   событий: давнего и очень давнего? 

Вопрос открыт --  над ним могут размышлять  философы,  а математикам   при некоторых   дополнительных  условиях   удается получить ответ.

%\newpage
\bf Результаты. \rm Кратно перемешивают   перемешивающие

гауссовские системы \cite{Leo};

действия, обладающие эргодической гомоклинической группой \cite{19}; 

автоморфизмы  с сингулярным спектром \cite{Ho};

унипотентные потоки \cite{Ma}--\cite{St}; %\cite{91},\cite{Mo}

потоки, обладающие   свойством Ратнер \cite{Ra} или его  аналогами \cite{FK}; 

квази-простые потоки \cite{RT};

 потоки положительного локального ранга \cite{00}; 

автоморфизмы конечного ранга  \cite{Ka}--  \cite{93}. 

\vspace{2mm}
Кратное перемешивание   устанавливается    путем нахождения свойств (алгебраических, спектральных, аппроксимационных),  вынуждающих  перемешивающую систему перемешивать кратно.  Наша цель  -- показать   связь  кратного перемешивания с некоторыми  такими свойствами. 

\%newpage

\section{Гауссовские и пуссоновские системы. Гомоклиническая группа} 
В теории представлений и в эргодической теории изучаются гауссовские и пуассоновские действия (см., например, книги \cite{KSF} и \cite{N}).  Первые индуцируются естественным вложением ортогональной группы в группу автоморфизмов  пространства с гауссовской мерой, а вторые являются непрерывным инъективным образом группы автоморфизмов  пространства c сигма-конечной мерой  в группе автоморфизмов пространства конфигураций, наделенного  пуассоновской мерой.
Для гауссовских перемешивающих автоморфизмов  кратное перемешивание доказал В.П. Леонов \cite{Leo}, а для пуассоновских перемешивающих систем оно непосредственно вытекает из определений и следующего  факта. Свойство "перемешивания" автоморфизма $T$  пространства с бесконечной мерой означает асимптотическую дизъюнктность  множеств $B_0$ и  $T^mB_1$  при $m\to\infty$ (множества $B_0,B_1$ имеют конечную меру).   Тогда при
$m_1,\dots m_k\, \to\,+\infty$   множества 
$B_0,  T^{m_1}B_1,  \dots , T^{m_1+\dots+m_k}B_k$ асимптотически дизъюнктны.  Пуассоновское соответствие превращает почти дизъюнктность этих множеств  в почти независимость отвечающих им цилиндрических множеств и тем самым обеспечивает настоящее кратное перемешивание для пуассоновской надстройки $P(T)$.  Можно  сказать, что гауссовское соответствие делает примерно то же самое, только почти дизъюнктность надо заменить на почти ортогональность, которую гауссовское соответствие превратит в асимптотическую  независимоcть. 

Мы  не будем  давать определения    гауссовских и пуассоновских систем, но скажем, что  они обладают эргодической  группой Гордина \cite{19},  которая обеспечивает свойство кратного перемешивания системы. Ограничимся обсуждением  этого факта. 

Пусть $T$ -- автоморфизм вероятностного пространства,  гомоклиническая группа Гордина $H(T)$ определяется так:
$$H(T)=\{S: \ T^{-n}ST^n\to I,\ n\to\infty\}.$$ 
М.И. Гордин показал, что эргодичность влечет за собой свойство перемешивания, а автор заметил, что  обычное перемешивание можно заменить на  кратное перемешивание.

\vspace{2mm}
\bf  Теорема 2.1. \it  Автоморфизм с эргодической гомоклинической  группой перемешивает с бесконечной кратностью. 
  \rm

\vspace{2mm}
Доказательство. Пусть $S\in H(T)$ -- эргодический автоморфизм.  Фиксируем целое $i$. Для любых измеримых множеств $A,B,C$  замечаем, что 
$$
\mu(S^{-i}A\cap T^{m} B\cap T^{m+n}C)-  
\mu(A\cap T^{m} B\cap T^{m+n}C)\to 0,\ m, n\to +\infty,$$
что следует из сходимости
$$\mu(A\cap T^{m}T^{-m} S^iT^{m} B\cap T^{m+n}T^{-m-n}S^iT^{m+n}C) -
\mu(A\cap T^{m} B\cap T^{m+n}C)\to 0,$$
которую в силу гомоклиничности автоморфизма $S$ обеспечивают условия  $$\mu (T^{-m} S^iT^{m}B\Delta B)\to 0, \ \ \ T^{-m-n}S^iT^{m+n}C
\Delta C)\to 0.$$ Но тогда выполняется  
$$
 \mu(A\cap T^{m} B\cap T^{m+n}C) -\frac 1 N \sum_{i=1}^N S^p
\mu(S^{-i}A\cap T^{m} B\cap T^{m+n}C)\to 0,\  m, n\to +\infty.$$
Теперь замечаем, что  для всякого $\eps>0$ найдется такое  $N$, что  для всех достаточно больших $ m,n$  выполнено 
$$
 \left|\mu(A)\mu(B)\mu(C) - \frac 1 N \sum_{i=1}^N 
\mu(S^{-i}A\cap T^{m} B\cap T^{m+n}C)\right| < \eps.$$
Действительно, 
из эргодичности автоморфизма $S$  имеем
$$ \left\|\frac 1 N \sum_{i=1}^N S^{-i}\chi_A - \mu(A)\chi_X\right\| \to \, 0, N\to\infty,$$
а из  свойства перемешивания автоморфизма  $T$ вытекает 
$ \mu(T^{m} B\cap T^{m+n}C)\to \mu(B)\mu(C).$
Таким образом, получаем
$$\mu(A\cap T^m B\cap T^{m+n}C)\to \mu(A)\mu(B)\mu(C), \ m,n\to +\infty.$$
Мы установили свойство перемешивания кратности 2.
Аналогично с использованием  индукции доказывается, что $T$ обладает  перемешиванием кратности $k>2$.

\bf Замечание.  \rm Группа $H(T)$  является полной (full). В силу теоремы Дая, если такая группа эргодична, она  содержит с точностью до сопряжения все эргодические автоморфизмы  (см. \cite{19}). 
Было бы интересно найти  неизоморфные  гауссовским и пуассоновским надстройкам  примеры систем с нулевой энтропией и эргодической  гомоклинической группой.  Вероятно, к ним относится давний  пример Дж. Кинга \cite{KG}.

\section{Самоприсоединения, марковские сплетения и  спектр. Теорема Оста.}
Самоприсоединения (self-joinings), а точнее говоря,
некоторые свойства действий, формулируемые в терминах  джойнингов или
марковских сплетающих операторов, являются одним из основных иструментов доказательства свойства кратного перемешивания. Приведем необходимые определения. 

Присоединением (или джойнингом) набора действий $\Psi_1,\dots,\Psi_n$ называется мера на $X^n =X_1\times \dots\times  X_n$  ($X_i=X$), проекции которой на ребра  куба $X^n$ равны $\mu$, причем мера инвариантна относительно 
диагонального действия произведения $\Psi_1\times \dots\times  \Psi_n$. Если $\Psi_1,\dots,\Psi_n$ суть копии одного  действия, такой джойнинг называется самоприсоединением.

Говорим, что действие $\Psi$  принадлежит классу $S(m, n)$, $n>m>1$  (или обладает свойством $S(m,n)$), если всякое cамоприсоединение  порядка $n>2$ такое, что все проекции на  $m$-мерные грани  куба $X^n$ равны $\mu^m$,  является тривиальным,  т.е. совпадает с мерой $\mu^n$,
произведением $n$ копий меры $\mu$.  Далее используем обозначение   $S_n=S(n-1,n)$, $n>2$.

\bf Джойнинги и кратное перемешивание. \rm Полезное наблюдение, к которому пришли независимо несколько исследователей, связывает свойства $S_n$  с кратным перемешиванием:

\it если перемешивающий автоморфизм обладает перемешиванием кратности 
$n-1$ и свойством $S_{n+1}$, то он обладает кратным перемешиванием
порядка $n$. \rm

Пусть для всяких измеримых $A,B,C$ выполняется
$$\mu(A\cap T^{m_i} B\cap T^{m_i+n_i}C)\to \nu(A\times B\times C), $$
Свойство $S_3$  автоморфизма $T$ влечет за собой
$$\nu(A\times B\times C) =\mu(A)\mu(B)\mu(C).$$
Действвительно,  $\nu$ --  нормированная мера на полукольце цилиндров 
вида $A\times B\times C$.  
Проверяется, что  
$$\nu(A\times B\times C) = \nu(TA\times TB\times TC),$$ 
$$\nu(A\times B\times X)=\mu(A)\mu(B), \ \nu(X\times B\times C) = \mu(B)\mu(C),$$
$$ \nu(A\times X\times C) = \mu(A)\mu(C).$$
Такие меры называют самоприсоединениями (джойнингами) с попарной независимостью.
Мере  $\nu$  однозначно соответствуе марковский оператор $$P_2::L_2^{\otimes 2}\to  L_2, \ \
(A, P_2(B\otimes C))=\nu(A\times B\times C)$$
(левые $A,B,C$ обозначают индикаторы  множеств $A,B,C$).
Из инвариантности  меры $\nu$ и ее проекционных свойств вытекает условие сплетения
$$ TP_2=P_2(T\otimes T), $$  и включение
$$ P^\ast_2  H \subset H\otimes H,$$
где  $H$ -- пространство функций,  ортогональных константам.

Если перемешивающий автоморфизм $T$ обладает свойством $S_3$, то для всяких измеримых множеств

$\nu(A\times B\times C)=\mu(A)\mu(B)\mu(C)$, тогда    $P^\ast_2  H=\{0\}$,
значит, и мы получаем, что $T$ обладает перемешивание кратности 2.  Аналогично получаем более общий факт:
Если перемешивающий автоморфизм $T$ обладает свойствами $S_3,\dots,S_{n+1}$, $n>2$,
то $T$ обладает перемешивание кратности $n$.

Приятным сюрпризом для специалистов было  обнаружение    Дж. Кингом   \cite{Ki} того, что одновременное выполнение свойств $S_3$ и $S_4$ влечет за собой выполнение свойств $S_p$ для всех $p>4$.
В \cite{97S}  показано, что   свойства $S_{2m}$, $m>1$, эквивалентны между собой и  вынуждают все свойства  $S_{2m-1}$, $m>1$.  Осталась неизученной связь перемешивания кратности $k$ с нечетными  свойствами $S_{2m-1}$. Впрочем, вопрос о существовании  автоморфизма со свойством $S_{2m-1}$, но без свойства $S_4$ следует отнести к трудным вопросам теории джойнингов.  

\bf Проблема дель Джунко-Рудольфа \cite{JR}: \it существует ли перемешивающий автоморфизм с нулевой энтропией, не обладающий свойством $S_p$? \rm  Этот вопрос также открыт.

Положительная энтропия и наличие дискретного спектра несовместимы со свойствами $S_p$. В работах \cite{RT}, \cite{97},\cite{GHR} для некоторых  классов действий показано, что свойство $S_3$    эквивалентно свойству $S_4$.
Известны  групповые действия, обладающие  свойствами   $S_{2m-1}$, но не обладают свойством  $S_4$. Пример: рассмотрим компактную коммутативную группу $X=\Z_2^{\Z}$ с мерой Хаара $\mu$, на $(X,\mu)$  определим действие $\Psi$, порожденное всеми сдвигами на группе $X$ и всеми автоморфизмами группы $X$. Оказывается, действие $\Psi$  и многие его конечно порожденные некоммутативные поддействия обладают всеми свойствами  $S_{2m-1}$, но не обладают свойством  $S_{4}$.

\vspace{2mm}
\bf Теорема (Ост).
\it Если автоморфизм $T$ обладает непрерывным сингулярным  спектром, то он обладает всеми свойствами $S_p$.  \rm

\vspace{2mm}
\bf Следствие.  \it Перемешивающий автоморфизм $T$ с сингулярным  спектром обладает перемешиванием всех кратностей. \rm

Вкратце схема рассуждения Оста следующая (в работе \cite{Ho}  рассматривается более общая ситуация джойнинга трех преобразований $R,S,T$, но мы ограничимся частным случаем, когда 
 $R=S=T$). 
Рассмотрим меру $\tau$ на торе $\T^3$:
$$\hat{\tau}(n_1,n_2,n_3) =\left(P_2(T^{n_1}f_1\otimes T^{n_2}f_2)\,,\, T^{n_3}f_3\right),$$
для фиксированных функций  $f_1,f_2,f_3$ с нулевым средним, а $P_2$ -- марковский оператор
(условного математического ожидания), отвечающий самоприсоединению $\nu$.

Носитель такой меры $\tau$ лежит в подгруппе $\{(t_1,t_2,t_3): t_1+t_2+t_3=0\}$,
причем  $\pi_{i,j}\tau\ll \sigma\otimes\sigma$, где   $\pi_{i,j}$ -- проекция  тора $\T_{1}\times\T_{2}\times \T_{3}$  на грань $\T_{i}\times \T_{j}$, $1\leq i< j\leq 3$, а $\sigma$ -- мера максимального спектрального типа автоморфизма $T$.

Методами гармонического анализа доказывается (и это ключевое место), что проекция $\pi\tau$ такой меры $\tau$   на $\T_1$ являтся суммой   дискретной  и абсолютной непрерывной мер,  при этом $\pi\tau\ll \sigma$.  Но $\sigma $ -- непрерывная сингулярная мера, значит, 
$\pi\tau=0$ , $\tau=0$.  Получаем, что для всех наборов функций с нулевым средним $f_1,f_2,f_3$ имеет место равенство
$$(P_2(T^{n_1}f_1\otimes T^{n_2}f_2)\, , \, T^{n_3}f_3)=0.$$ 
Следовательно, 
$$ P_2  (H \otimes H)  =\{0\}.$$
Таким образом, самоприсоединение $\nu$, с которым связан оператор $P_2$, тривиально,
тем самым доказано свойство $S_3$.

Ост доказал больше: для всякого джойнига автоморфизмов $R,S,T$ 
с попарной независимостью, из сингулярности спектра $T$ вытекает тривиальность джойнинга.
Отсюда следует, что $T$ обладает всеми свойствами $S_p$ (для доказательства $S_4$ нужно 
рассмотреть случай  $R=T\times T$, $S=T$).   

\section{Взаимосвязи свойств $\bf S_p$.} \rm Ниже  мы рассмотрим    частный случай  теоремы Оста и  покажем, как используются  джойниги и марковские сплетения для выявления   взаимодействия свойств $S_p$.

\vspace{2mm}
\bf  Утверждение   4.1. \it 
 Если спектральная мера $\sigma$ перемешивающего автоморфизма $T$  не имеет общей компоненты со сверточной степенью $\sigma^{\ast 5}$,  то автоморфизм $T$  обладает всеми свойствами $S_p$ и по этой причине обладает перемешиванием всех кратностей.  \rm

\vspace{2mm}
Следует сказать, что такая мера $\sigma$ обязана быть сингулярной,
поэтому результат следует из теоремы Оста.  Наша задача -- в частном случае показать простой метод, использующий  алгебру джойнингов. 
Докажем  перемешивание кратности 2. Пусть для всяких измеримых $A,B,C$ выполняется
$$\mu(A\cap T^{m_i} B\cap T^{m_i+n_i}C)\to \nu(A\times B\times C), $$
наша цель -- показать, что
$$\nu(A\times B\times C) =\mu(A)\mu(B)\mu(C).$$
Проверяетсяется, что  $\nu$ --  мера, удовлетворяющая  свойствам:
$$\nu(A\times B\times C) = \nu(TA\times TB\times TC),$$ 
$$\nu(A\times B\times X)=\mu(A)\mu(B), \ \nu(X\times B\times C) = \mu(B)\mu(C),$$
$$ \nu(A\times X\times C) = \mu(A)\mu(C),$$
Такие меры называют самоприсоединениями (джойнингами) с попарной независимостью.
Мере  $\nu$  однозначно соответствуе марковский оператор $$P_2::L_2^{\otimes 2}\to  L_2, \ \
(A, P_2(B\otimes C))=\nu(A\times B\times C)$$
(левые $A,B,C$ обозначают индикаторы правых множеств $A,B,C$).
Из инвариантности  меры $\nu$ и ее проекционных свойств вытекает условие сплетения
$$ TP_2=P_2(T\otimes T), $$  и включение
$$ P^\ast_2  H \subset H\otimes H,$$
где  $H$ -- пространство функций,  ортогональных константам.

Рассмотрим операторы $P_3:L_2^{\otimes 3}\to  L_2$ и $P_5:L_2^{\otimes 5}\to  L_2$, которые заданы следующим образом:
 $$( P_3(A_1\otimes A_2\otimes A_3)\,,\, A_4)=(P_2(A_1\otimes A_2)\,,\,  P_2(A_3\otimes A_4)),$$
$$( P_5(A_1\otimes \dots\otimes A_5)\,,\, A_6)=
(P_3(A_1\otimes A_2\otimes A_3),  P_3(A_4\otimes A_5\otimes A_6)).$$
Проверяется, что имеет место включение 
$$P_5^\ast H\subset H^{\otimes 5}$$
и выполнено условие сплетения
$$   PT^{\otimes 5}=TP.$$
Из $P_5^\ast H\neq {0}$ следует, что спектральная мера $\sigma$ имеет общую компоненту 
со сверткой $\sigma^{\ast 5}$.  Но по условиям теоремы эти меры  дизъюнктны, значит,
 $$P_5H^{\otimes 5}={0}.$$
Это влечет за собой $$P^\ast_3P_3H^{\otimes 3}=0, \ \  P_3H^{\otimes 3}={0},$$ 
тогда 
$$P^\ast_2 P_2H^{\otimes 2}={0}, \ \  P_2H^{\otimes 2}={0},$$ 
значит,
$$(A, P_2(B\otimes C))= \mu(A)\mu(B)\mu(C)=\nu(A\times B\times C).$$
Тем самым мы доказали, что автоморфизм $T$ обладает перемешиванием кратности 2.
Аналогичным методом  доказываем перемешивание кратности 3, 4 и т.д.,
и даже более общее утверждение: \it 
 если для перемешивающего  автоморфизма $T$  со  спектральной мерой $\sigma$ 
 некоторые  сверточные степени $\sigma^{\ast m}$ и  $\sigma^{\ast n}$ взаимно
сингулярны, то $T$   перемешивает с бесконечной кратностью. \rm

\vspace{2mm}
   Следующая теорема обобщает результат  Дж. Кинга из \cite{Ki}.

\vspace{2mm}
\bf Теорема 4.2 ( \cite{97S}). \it Если для сохраняющего меру действия  найдется нетривиальное самоприсоединение порядка $n>2$
с попарной независимостью, то  для него и  всякого $p>1$ найдется нетривиальное самоприсоединение класса  $M(2p-1,2p)$.  \rm

\vspace{2mm} \bf
Следствие. \it  Пусть  перемешивающее действие коммутативной группы для некоторого $k>1$ не
обладает перемешиванием кратности  $k$, тогда для всякого $p>1$ оно допускает нетривиальное самоприсоединение класса  $M(2p-1,2p)$.\rm

\vspace{2mm}
 Теорема 3.2 доказывается тем же методом, что и  теорема 3.1.
Из нетривиальной меры $\nu_2=\nu$, самоприсоединения класса $M(2,3)$, мы
изготавливаем нетривиальное  самоприсоединение $\nu_5$ класса $M(5,6)$:
$$\nu_5(A_1\times \dots\times A_6):=
(P_3(A_1\otimes A_2\otimes A_3),  P_3(A_4\otimes A_5\otimes A_6)).$$
При этом  порядок самоприсоединения возрос. Аналогичным образом можно понижать порядок (но не ниже порядка 4).  Действительно,  пусть самоприсоединение $\nu$ имеет порядок $p+2$, причем проекции меры $\nu$ на все $p+1$-мерны грани суть $\mu^{p+1}$ (мера класса $M(p+1,p+2)$.  Определим оператор $P:L_2^{\otimes 2}\to  L_2^{\otimes p}$ формулой 
$$(P(A_1\otimes A_2), B_1\otimes \dots \otimes B_p)=\nu(A_1\times A_2\times 
B_1\times \dots \times B_p).$$
Положим
$$\nu_2(A_1\times A_2 \times A_1'\times A_2'):= (P(A_1\otimes A_2), P(A_1'\otimes A_2')).$$
Если $\nu_2= \mu^{\otimes 4}$, тогда получим  $\nu= \mu^{\otimes(p+2) }$ по тем же причинам, что и в доказательстве утверждения 3.1.

\section{Коммутационные соотношения}
Если автоморфизм обладает лебеговским спектром и нулевой энтропией
(например, является элементом  орициклического или унипотентного потока), доказательство  свойства кратного перемешивания  требует других подходов.  
Если автоморфизм  является элементом действия некоммутативной группы Ли,  как правило,  коммутационные соотношения вынуждают   свойство кратного перемешивания.  В качестве интересного примера рассмотрим  действие группы Гейзенберга.
\medskip

 \bf Теорема 5.1. \it Пусть сохраняющие меру эргодические потоки
$\Psi_r ,\Phi_s, T_t $  удовлетворяют соотношению
$$\Psi_r\Phi_s= T_{rs}\Phi_s\Psi_r,$$
причем   $T_c$ -- слабо перемешивающий поток, коммутирующий с
 $\Psi_r$ и $\Phi_s$. Тогда потоки $\Psi_r$ и $\Phi_s$ обладают перемешиванием  всех порядков. \rm
\medskip

 Доказательство. Определим  оператор $P$:
$$
  \left< A,P(B\otimes C)\right>  =
\lim_{i'\to\infty} \left< A,\,\Psi_{m(i')}B \Psi_{n(i')}C)\right>.
$$
Предположим, что $ 0 < m(i') < n(i')$ и существует предел
 $$\lim_{i} m(i)/n(i) = \alpha$$
для $i$ --  подпоследовательности
 $i'$.

Заметим, что при $s\to 0$ выполнено
$$
\left< \Phi_sA,\,(\Psi_{m(i)}\Phi_sB)(\Psi_{n(i)}\Phi_sC)\right>  \to
\left< A,\,(\Psi_{m(i)}B)(\Psi_{n(i)}C)\right>.
$$
Следовательно,
$$
\lim_{\eps\to 0}\lim_{i\to\infty} \frac{1}{\eps}\int_{0}^\varepsilon
\left< \Phi_sA,\,(\Psi_{m(i)}\Phi_sB)(\Psi_{n(i)}\Phi_sC)\right> ds =
$$
$$
\lim_{i\to\infty} \left< A,\,(\Psi_{m(i)}B)(\Psi_{n(i)}C)\right>.
$$
Ввиду коммутационных соотношений  получим
$$
\left< \Phi_sA,\,(\Psi_{m(i)}\Phi_sB)(\Psi_{n(i)}\Phi_sC)\right>  =
\left< A,\,(\Psi_{m(i)}T_{sm(i)}B)(\Psi_{n(i)}T_{sn(i)}C)\right>.
$$
Таким образом,
$
\left< A,P(B\otimes C)\right> =
$
$$
\lim_{\eps\to 0}\lim_{i\to\infty} \frac{1}{\eps}\int_{0}^\varepsilon
\left< A,\Psi_{m(i)}T_{sm(i)}B
\Psi_{n(i)}T_{sn(i)}C)\right>  ds.
$$
 Выражение
$$ \frac{1}{\eps}\int_{0}^\varepsilon
\left< A,\Psi_{m(i)}T_{sm(i)}B
\Psi_{n(i)}T_{sn(i)}C)\right>  ds
$$
при фиксированном $r$ и больших $i$ мало отличается от
$$ \frac{1}{\eps}\int_{0}^\eps
\left< A,\Psi_{m(i)}T_{(s-\frac{r}{n(i)})m(i)}B
\Psi_{n(i)}T_{(s-\frac{r}{n(i)})n(i)}C)\right>  ds.
$$
Поэтому при фиксированном произвольном $r\in\bf R$ получим
$
\left< A,P(B\otimes C)\right> =
$
$$
\lim_{\eps\to 0}\lim_{i\to\infty} \frac{1}{\eps}\int_{0}^\varepsilon
\left< A,\Psi_{m(i)}T_{sm(i)}T_{ar}B
\Psi_{n(i)}T_{sn(i)}T_{r}C)\right>  ds.
$$
Это  приводит к  равенству
$$
P(T_{ar}\otimes T_{r} )= P,
$$
что влечет за собой
$$\left< A,P(B\otimes C)\right> = \mu(A)\mu(B)\mu(C).
$$

\vspace{2mm}
\it Верно ли что действие центра группы Гейзенберга 
(т.е. поток $T_{t}$, коммутирующий с потоками $\Psi_r, \Phi_s$) обладает кратным перемешиванием?\rm

Кратное перемешивание некоторых потоков как следствие  нетривиальных коммутационных соотношений было получено в работах \cite{Ma}-- \cite{Mo}.

\section{Количество отклонений от кратного перемешивания. Локальный ранг.}

\bf $\bf 1^+$-перемешивание. \rm  Определим   свойство
$Mix(1^+)$, занимающее промежуточное положение между свойством
перемешивания кратности 1 и кратности 2.
Положим
$$Der(\varepsilon ,A,B,C,h) = 
\lbrace (z,w)\in Q(\varepsilon ,h): | \mu (A\cap T^z B\cap
T^w C) - \mu (A)\mu (B)\mu (C)| > \varepsilon \rbrace ,
$$
где
$$
  Q(\varepsilon ,h) = \lbrace (z,w) \epsilon [0,h]^2 :
|z|,|w|,|z-w| > \varepsilon h \rbrace ,
$$
$$
  dev(h)=  |\,Der(\varepsilon, A,B,C, h)\,|/ h.
$$

 Если  автоморфизм $T$ перемешивает двукратно, то $dev(h)=0$
для больших значений $h$, \it если $T$ перемешивает, то последовательность $dev(h)$ ограничена. \rm Это легкое следствие (см. \S 2, \cite{92}) известной теоремы Блюма-Хансона о том, что
для перемешивающих систем средние вдоль возрастающих последовательностей сходятся по норме к константе. 

\it Пишем  $T\in Mix(1^+)$,  если $dev(h) \to \, 0$, $h\to\infty$, 
для любых измеримых множеств $A,B,C$ и $\varepsilon >0.$   \rm

 \bf Локальный ранг автоморфизмов. \rm
Говорят, что автоморфизм $T$ обладает \it локальным рангом \rm
$\beta=\beta(T) > 0$, если
для некоторой последовательности $U_j=\bigsqcup_{k\in Q_j}T^zB_j$ башен
Рохлина-Халмоша, где $Q_j=\{0,1,\dots,h_j\}$, выполнены условия:
$\mu(U_j)\to \beta$ и для каждого измеримого $A\in\B$ пересечение  $U_j\cap A$ асимптотически как угодно точно приближается объединением некоторых  этажей в башне. 

С.Каликов \cite{Ka} доказал, что перемешивающие автоморфизмы $T$ в случае $\beta(T)=1$ обладают перемешиванием кратности 2. 
Теорема Каликова только частично была обобщена на автоморфизмы с положительным локальным рангом.

\vspace{2mm}
\bf Теорема 6 (\cite{95}). \it Пусть $T\in Mix(1^+)$ и $\beta(T)>0$.
 Тогда $T\in Mix(\infty)$.   \rm

\vspace{2mm}
\bf Гипотетический пример \rm  перемешивающего автоморфизма $T\in Mix(1)\setminus Mix(2)$ положительного локального ранга обладает весьма интересными свойствами.
Для него $dev(h)\sim h$ (отклонения "максимально" часты), найдутся множества $A,B$ положительной меры и последовательности $m_i,n_i\to +\infty$ такие, что
$$A\cap T^{m_i}A\cap T^{m_i+n_i}B =\varnothing,$$ \cite{97},
а спектр автоморфизма лебеговский  конечной кратности (вытекает из \cite{GHR}). Не думаю,  что существует  такой автоморфизм, а в случае  $\beta(T)>1/2$ его точно нет \cite{00}.

\vspace{4mm}
\bf Благодарности. \rm Автор  признателен  участникам  семинара по многомерному комплексному анализу за полезные обсуждения тематики "алгебра джойнингов".

\end{fulltext}
%\newpage

\end{document}